\documentclass[MSNbibl,number,citesort,seceqn,dvips]{arxbj}
\usepackage{upgreek}

\aid{0}
\volume{19}
\issue{3}
\pubyear{2013}
\firstpage{1006}
\lastpage{1027}
\doi{10.3150/12-BEJ415} 

\makeatletter
\newtheorem{them}{Theorem}
\newtheorem{proposition}{Proposition}
\newtheorem{lemma}{Lemma}
\newremark{remark}{Remark}
\newcommand{\ep}{\epsilon}

\renewcommand{\)}{)}
\newcommand{\no}{\nonumber}
\newcommand{\la}{\label}
\newcommand{\be}{\begin{eqnarray}}
\newcommand{\ee}{\end{eqnarray}}
\newcommand{\bestar}{\begin{eqnarray*}}
\newcommand{\eestar}{\end{eqnarray*}}
\newcommand{\al}{\alpha}
\newcommand{\pr}{\mathsf{P}} 
\renewcommand{\ep}{\mathsf{E}} 
\def\Var{\mathsf{Var}} 
\renewcommand{\limsup}{\overline{\lim}}
\renewcommand{\liminf}{\underline{\lim}}
\newcommand{\F}{\mathbf{F}}
\newcommand{\E}{\mathbf{E}}
\newcommand{\A}{\mathbf{A}}
\newcommand{\C}{\mathbf{C}}
\newcommand{\nn}{\nonumber}
\newcommand{\lbl}{\label}
\newcommand{\eq}[1]{(\ref{#1})}

\makeatother

\begin{document}
\begin{frontmatter}

\title{Self-normalized Cram\'{e}r type moderate deviations for the
maximum of sums}
\runtitle{Self-normalized moderate deviation for maximum}

\begin{aug}
\author[1]{\fnms{Weidong} \snm{Liu}\thanksref{1}\ead[label=e1]{liuweidong99@gmail.com}},
\author[2]{\fnms{Qi-Man} \snm{Shao}\corref{}\thanksref{2}\ead[label=e2]{maqmshao@ust.hk}}
\and
\author[3]{\fnms{Qiying} \snm{Wang}\thanksref{3}\ead[label=e3]{qiying@maths.usyd.edu.au}}
\runauthor{W. Liu, Q.-M. Shao and Q. Wang} 
\address[1]{Department of Mathematics and Institute of Natural
Sciences, Shanghai Jiao Tong University, Shanghai, China. \printead{e1}}
\address[2]{Department of Mathematics, Hong Kong University of Science
and Technology, Clear Water Bay, Kowloon, Hong Kong. \printead{e2}}
\address[3]{School of Mathematics and Statistics, University of Sydney,
Australia.\\ \printead{e3}}
\end{aug}

\received{\smonth{10} \syear{2010}}
\revised{\smonth{10} \syear{2011}}

%
\begin{abstract}
Let $X_{1}, X_{2},\ldots$ be independent random variables with
zero means and finite variances, and let $S_{n}=\sum_{i=1}^{n}X_{i}$
and $V^{2}_{n}=\sum_{i=1}^{n}X^{2}_{i}$. A Cram\'{e}r type moderate deviation
for the maximum of
the self-normalized sums $\max_{ 1 \leq k \leq n} S_k/V_n$
is obtained. In particular, for identically distributed $X_1,
X_2,\ldots
,$ it is proved that
$\textsf{P}(\max_{1\leq k\leq n}S_{k}\geq x V_{n})/(1-\Phi
(x))\rightarrow2$
uniformly for $0<x\leq\mathrm{o}(n^{1/6})$ under the optimal finite
third moment of $X_1$.
\end{abstract}

%

\end{frontmatter}

\section{Introduction and main results}
\label{intro.sec}
Let $X_1, X_2,\ldots$ be a sequence of independent non-degenerate random
variables with zero means. Set
\[
S_n=\sum_{j=1}^nX_j \quad\mbox{and}\quad V_n^2=\sum_{j=1}^nX_j^2.
\]
The past decade has brought significant developments
in the limit theorems for the so-called ``self-normalized'' sum, $S_n/V_n$.
It is now well understood that the limit theorems for $S_n/V_n$ usually
require fewer moment assumptions than those for
their classical standardized counterpart, and thus have much wider
applicability. For examples, for identically distributed $X_1,
X_2,\ldots,$
a self-normalized large deviation holds without any moment assumption
 (Shao \cite{Shao97}),
and a Cram\'er type moderate deviation (Shao \cite{Shao99}),
\begin{equation}\label{in1}
\lim_{n \rightarrow\infty}
\frac{ \pr( S_n \ge x V_n )}{1-\Phi(x)}=1,
\end{equation}
holds uniformly for $x \in[0, \mathrm{o}(n^{1/6}))$ provided that
$\ep|X_1|^3 < \infty$,
whereas a finite moment-generating condition of $\sqrt{|X_1|}$
is necessary for a similar result for the
standard sum $S_n/\sqrt{\Var( S_n)}$ (see, e.g., Linnik \cite{Lin62}). For
more related results, we
refer to de la Pe\~{n}a, Lai and Shao \cite{PLS09} for a systematic
treatment of the theory and applications of self-normalization
and Wang \cite{W2} for some refined self-normalized moderate deviations.%

As for the Cram\'{e}r type moderate deviations for the maximum of
self-normalized sums, namely for $\max_{1\le k\le n}S_k/V_n$,
Hu, Shao and Wang \cite{HSW09} were the first to prove that if $X_1,
X_2,\ldots$ is a sequence
of i.i.d. random variables with $\mathsf{E}X_1^4<\infty$, then
\begin{equation}\label{t1.1a}
\lim_{n \rightarrow\infty} \frac{ \pr(\max_{1\le k \le n} S_k
\ge x V_n )}{1-\Phi(x)}=2,
\end{equation}
uniformly for $x \in
[0, \mathrm{o}(n^{1/6}))$.
This contrasts with the moderate deviation result for the maximum of
partial sums of Aleshkyavichene
\cite{Ale79a,Ale79b}, where a finite moment-generating condition is required.
However, in view of the result given in (\ref{in1}),
it is natural to ask whether a finite third moment suffices for (\ref{t1.1a}).
The main purpose of this paper is to provide an affirmative answer to
this question. Indeed, we have the following more general result for
independent random variables.
\begin{them} \label{th0}
Assume that $\max_{k\geq1}\ep|X_{k}|^{2+r}<\infty$ and $\min
_{k\geq
1}\ep X^{2}_{k}>0$, where $0<r\le1$.
Then (\ref{t1.1a}) holds uniformly in $0\le x\leq \mathrm{o}(n^{r/(4+2r)})$.
\end{them}

As in the moderate deviation result for self-normalized sum $S_n/V_n$,
Theorem~\ref{th0} is sharp in both the moment condition and the range
in which the result (\ref{t1.1a}) holds true. Examples can be
constructed similarly as done by Chistyakov and G\"otze \cite{CHGPG03} and
Shao~\cite{Shao99}. In particular, for $r=1$ and identically distributed $X_1,
X_2,\ldots$\,, Theorem~\ref{th0} establishes (\ref{t1.1a}) under the
optimal finite third moment of $X_1$.

Theorem~\ref{th0} can be extended further; in fact, it is a direct
consequence of Theorem~\ref{th1} below.
Set $B^{2}_{n}=\sum_{i=1}^{n}\ep X^{2}_{i}$, $L_{n,r}=\sum
_{i=1}^{n}\ep
|X_{i}|^{2+r}$ and $d_{n,r}=B_{n}/L^{1/(2+r)}_{n,r}$,
where $0< r\le1$.
\begin{them} \label{th1} For $0<r\leq1$,
suppose that $d_{n,r}\rightarrow\infty$ as $n\to\infty$, and that
\begin{equation}\label{ad2}
\max_{1\leq k\leq n}\frac{\sum_{j=k}^{n}\ep
|X_{j}|^{2+r}}{\sum_{j=k}^{n}\ep|X_{j}|^{2}} \leq \frac{\tau
L^{r/(2+r)}_{n,r}}
{d^{\delta}_{n,r}}\qquad
\mbox{for some } \delta, \tau>0.
\end{equation}
Then (\ref{t1.1a}) holds uniformly in $0\leq x\leq\min\{B_n,
\mathrm{o}(d_{n,r})\}$.
\end{them}
\begin{remark}
For i.i.d. random variables with $\ep X_i =0$ and $\ep|X_i|^3 < \infty$, Jing,
Shao and Wang \cite{JSW03}
proved that \eq{in1} can be refined as
\[
{ \pr(S_n \geq x V_n) \over1- \Phi(x)} = 1 + \mathrm{O}(1) (1+x^3) \ep|X_1|^3 /
(\ep X_1^2)^{3/2}
\]
uniformly in $x \in[0, n^{1/6} (\ep X_1^2)^{1/3} / (\ep X_1^2)^{1/2})$,
where $\mathrm{O}(1)$ is bounded by an absolute
constant. We conjecture that a similar result holds for $\max_{1 \leq k
\leq n} S_k /V_n$, that is,
\[
{\pr(\max_{ 1\leq k \leq n} S_k \geq x V_n) \over1- \Phi(x)} = 2 + \mathrm{O}(1)
(1+x^3) \ep|X_1|^3 / (\ep X_1^2)^{3/2}
\]
uniformly in $x \in[0, n^{1/6} (\ep X_1^2)^{1/2} / (\ep|X_1|^3)^{1/3})$.
\end{remark}

This paper is organized as follows. The proof of the main theorems is
given in the next section.
The proofs of two technical propositions are deferred to Sections \ref
{sec4} and~\ref{sec5}, respectively.
Throughout the paper, $A, A_1, \ldots$ denotes absolute constants
and $C_{\delta,\tau}$ denotes a constant depending only on $\delta$ and
$\tau$, which might be
different at each appearance.
%
\section{Proofs of theorems}
\begin{pf*}{Proof of Theorem~\ref{th0}} Simple calculations show that if
\[
\max_{k\geq1}\ep|X_{k}|^{2+r}<\infty \quad\mbox{and}\quad \min_{k\geq1}\ep
X^{2}_{k}>0,
\]
then $B_n^2\asymp n$, $L_{n,r}\asymp n$, $d_{n,r}\asymp
n^{r/(4+2r)}$ and (\ref{ad2}) holds for
$\delta=1$ and some $\tau>0$, where
the notation $a_n \asymp b_n$ denotes $0<\liminf_{n\to\infty}
a_n/b_n<\limsup_{n\to\infty}a_n/b_n<\infty$.
Therefore, Theorem~\ref{th0} follows immediately from Theorem~\ref{th1}.
\end{pf*}
\begin{pf*}{Proof of Theorem~\ref{th1}} First note that for $\forall
\epsilon>0$,
\[
\frac1{B_n^2}\sum_{k=1}^n\ep X_k^2I(|X_k|\ge\epsilon B_n) \le\epsilon
^{-r} d_{n,r}^{-1/(2+r)}\to0
\]
whenever $d_{n,r}\to0$. That is, the Lindeberg condition is satisfied
for the sequence $X_1, X_2,\ldots.$ On the other hand, routine
calculations show that, given $d_{n,r}\to0$,
$V_{n}^2/B_{n}^2\rightarrow1$ in probability. Given these facts, the
invariance principle (see Theorem 2 of Brown \cite{BrBM71}) and the continuous
mapping theorem imply that $\max_{1\le k\le n}S_k/V_n\to_D |N(0,1)|$.
This yields
(\ref{t1.1a}) uniformly for $0\leq x\leq M$, where $M$ is an arbitrary
constant.
Thus, Theorem~\ref{th1} will follow if we can prove
\be
\lim_{M\to\infty}\lim_{n\to\infty}\sup_{M\le x\leq\min\{B_n,
\mathrm{o}(d_{n,r})\}}
\biggl|\frac{\pr{(\max_{1\leq k\leq n}S_{k}\geq xV_{n} )}}{1-\Phi(x)}-2 \biggr|
=0. \la{ad190}
\ee
Toward this end, let
\[
\Delta_{n,x} = \frac{x^{2}}{B^{2}_{n}}\sum_{i=1}^{n}\ep
X^{2}_{i}\{|X_{i}|>B_{n}/x\}+\frac{x^{3}}{B^{3}_{n}}\sum_{i=1}^{n}\ep
|X_{i}|^{3}
I\{|X_{i}|\leq B_{n}/x\},
\]
and write
\begin{equation}\label{n0a}
n_{0}\equiv n_{0}(x)=\max\Biggl\{k\dvt \sum_{j=k}^{n}\ep X^{2}_{j}\geq192 B^{2}_{n}
\log(x\vee \mathrm{e})/x^2, 1\le k\le n \Biggr\}.
\end{equation}
It can be readily seen that the condition (\ref{ad2}), together with
$0<x\leq\min\{B_n,\mathrm{o}(d_{n,r})\}$ and $d_{n,r}\to\infty$, imply the
existence of an absolute constant $A$ such that
\begin{equation}\label{a0}
0\le x\le B_n,\qquad \Delta_{n,x}\leq\min(\delta^{9/2}, 1)/A,
\end{equation}
$\Delta_{n,x}\to0$, 
and
\begin{equation}\label{a2}
\frac{\sum_{j=n_{0}+1}^{n}\ep|X_{j}|^{3}I\{|X_{j}|\leq B_{n}/x\}
}{\sum
_{j=n_{0}+1}^{n}\ep|X_{j}|^{2}} \leq
\frac{B_{n}}{x^{1+\delta}}
\end{equation}
for all sufficiently large $n$, where $\delta$
is defined as in (\ref{ad2}). The result (\ref{ad190}) follows
immediately from the following proposition.
\begin{proposition} \label{th2} For all $x\ge2$ satisfying (\ref{a0})
and (\ref{a2}),
we have
\begin{equation}\label{ad3}
\frac{\pr{(\max_{1\leq k\leq n}S_{k}\geq xV_{n})}}{1-\Phi(x)}=
2+\mathrm{O}(1) \bigl(x^{-\min\{1/4, \delta/20\}}+\Delta_{n,x}^{1/9} \bigr),
\end{equation}
where
$\mathrm{O}(1)$ is bounded by a constant $C_{\delta}$ that depends only
on $\delta$.
\end{proposition}

The main idea of the proof of Proposition~\ref{th2}
is to use truncation and the maximum probability inequality and then apply
a moderate deviation theorem of Sakhanenko \cite{Sak92} to the truncated variables.
A suitable truncation level is ensured by using an inequality from
Jing, Shao and Wang \cite{JSW03}, page 2181. This avoids the conjugate argument
of Hu, Shao and Wang \cite{HSW09}, and makes it possible to prove the main
result under an optimal moment assumption.

It remains to prove Proposition~\ref{th2}. In addition to the notation
in the previous section, let
$\gamma=72^{-1}\min(\delta,1)$,
\[
\varepsilon= \max(2 \Delta^{2/9}_{n,x}, \gamma x^{-1/2}, \gamma
x^{-\delta/10}),\qquad m=[x^{2}/2],
\]
$N_0=\varnothing$ and, for $1\leq l\leq m$,
$N_l=\{j_{1},j_{2},\ldots,j_{l}\}\subseteq
\{1,2,\ldots,n\}$. Furthermore, write
$\bar{X}_{i}=X_{i}I\{|X_{i}|\leq\varepsilon B_{n}/x\}$, and for
$0\le l\le m$ and $ 1\le k\le n$,
\begin{eqnarray*}
\bar{S}_{k}^{N_l}&=&\sum_{i=1,i\notin N_l}^{k}\bar{X}_{i},\qquad
(\bar{V}^{N_l}_{n})^2=\sum_{i=1,i\notin N_l}^{n}\bar{X}^{2}_{i},\qquad
({\bar{B}}^{N_l}_{n})^2=\sum_{i=1,i\notin N_l}^{n}\ep{\bar
{X}}^{2}_{i},\\
 S_{k}^{N_l}&=&\sum_{i=1,i\notin N_l}^{k}X_{i},\qquad
(V^{N_l}_{n})^2=\sum_{i=1,i\notin N_l}^{n}X^{2}_{i},\qquad
({B}^{N_l}_{n})^2=\sum_{i=1,i\notin N_l}^{n}\ep{X}^{2}_{i}.
\end{eqnarray*}
Note that if $s,t\in R^{1}$, $x\geq1$, $c\geq0$ and $s+t\geq x\sqrt
{c+t^{2}}$, then $s\geq(x^{2}-1)^{1/2}\sqrt{c}$.
Similar to the arguments in reported by Jing, Shao and Wang \cite{JSW03},
page
2181, we have
\begin{eqnarray}\label{ad4}
\pr{\Bigl(\max_{1\leq k\leq n}S_{k}\geq
xV_{n}\Bigr )}
&\leq& \pr{\Bigl(\max_{1\leq k\leq n}\bar{S}_{k}^{N_0}\geq
x\bar{V}_{n}^{N_0} \Bigr)}\nonumber\\[-2pt]
&&{}+ \sum_{j_1=1}^{n}
\pr{\Bigl(\max_{1\leq k\leq n}
S_{k}\geq
xV_{n}, |X_{j_1}|\geq
\varepsilon{B_{n}}/x\Bigr )}\nonumber \\[-9pt]\\[-9pt]
&\leq& \pr\Bigl{(}\max_{1\leq k\leq
n}\bar{S}_{k}^{N_0}\geq x\bar{V}_{n}^{N_0} \Bigr{)} \nonumber\\[-2pt]
&&{} + \sum_{j_1=1}^{n}
\pr{\Bigl(\max_{1\leq k\leq n}S^{N_1}_{k} \geq
\sqrt{x^{2}-1} V^{N_1}_{n} \Bigr)}\pr{(|X_{j_1}|\geq
\varepsilon{B_{n}}/x )}\nonumber
\end{eqnarray}
and
\begin{eqnarray}\label{ad6}
\pr{\Bigl(\max_{1\leq k\leq n}S_{k}\geq
xV_{n}\Bigr )}
&\geq& \pr{\Bigl(\max_{1\leq k\leq
n}\bar{S}_{k}^{N_0}\geq x\bar{V}_{n}^{N_0}\Bigr )} \nn\\[-9pt]\\[-9pt]
&&{}-\sum_{j_1=1}^{n}\pr{\Bigl(\max_{1\leq k\leq n}\bar{S}^{N_1}_{k} \geq
\sqrt{x^{2}-1} \bar{V}^{N_1}_{n} \Bigr)}\pr{(}|X_{j_1}|\geq
\varepsilon{B_{n}}/x {)}.\nn
\end{eqnarray}
Repeating \eq{ad4} $m$-times gives
\begin{eqnarray}\label{eq5}
\pr{\Bigl(\max_{1\leq k\leq n}S_{k}\geq xV_{n} \Bigr)}
& \leq& \pr{\Bigl(\max_{1\leq k\leq n}\bar{S}_{k}^{N_0}\geq x\bar
{V}_{n}^{N_0} \Bigr)}\nn \\[-9pt]\\[-9pt]
& & {}+
\sum_{k=1}^{m}Z_{k}(x)+ {\Biggl\{\sum_{k=1}^{n}\pr{(}|X_{k}|\geq
\varepsilon B_{n}/x {)} \Biggr\}}^{m+1}, \nn
\end{eqnarray}
where
\[
Z_{k}(x)=\sum_{j_{1}=1}^{n}\cdots\sum_{j_{k}=1}^{n} {\Biggl[\prod
_{i=1}^{k}\pr{(}|X_{j_{i}}|\geq
\varepsilon B_{n}/x {)} \Biggr]} \times\pr{\Bigl(\max_{1\leq j\leq
n} {\bar S}^{N_k}_{j}\geq\sqrt{x^{2}-k} {\bar
V}^{N_k}_{n} \Bigr)}.
\]
Note that
\begin{eqnarray} \label{loser1}
&&\sum_{k=1}^{n}\pr{(}|X_{k}|\geq\varepsilon B_{n}/x {)} \nn\\[-2pt]
&&\quad \leq \frac{x^{2}}{\varepsilon^{2}B^{2}_{n}}\sum_{k=1}^{n}\ep
X^{2}_{k}I\{|X_{k}| \geq\varepsilon B_{n}/x\}\nn
\\[-9pt]
\\[-9pt]
&&\quad \leq \frac{x^{2}}{\varepsilon^{2}B^{2}_{n}}\sum_{k=1}^{n}\ep
X^{2}_{k}I\{|X_{k}|
\geq B_{n}/x\} + \frac{x^{3}}{\varepsilon^{3}B^{3}_{n}}
\sum_{k=1}^{n}\ep|X_{k}|^{3}I\{|X_{k}| \leq B_{n}/x\}\nn\\[-2pt]
&&\quad \leq \varepsilon^{-3}\Delta_{n,x} \leq\varepsilon^{3/2}/16\le
1/16.\nn
\end{eqnarray}
It follows from $m=[x^2/2]$ that
\begin{equation}
{\Biggl[\sum_{k=1}^{n}\pr{(}|X_{k}|\geq
\varepsilon B_{n}/x {)} \Biggr]}^{m+1}\leq \mathrm{e}^{-x^{2}}. \label{los0}
\end{equation}
This, together with (\ref{ad6}) and (\ref{eq5}), implies that
Proposition~\ref{th2} will follow if we prove the following two\vspace*{-1pt} propositions.
\begin{proposition}\label{lemma-1}
For all $0\le l\le m$, all $x/2\leq y\leq x$, and all $x\ge2$
satisfying (\ref{a0}) and~(\ref{a2}), we have
\begin{equation}
\frac{\pr{(}\max_{1\leq k\leq n}\bar{S}^{N_l}_{k}\geq
y\bar{V}^{N_l}_{n} {)}}{1-\Phi(y)}
\leq2+C_{\delta,\tau}(\varepsilon^{-2}\Delta_{n,x}+\varepsilon).
\label{ad8}
\end{equation}
\end{proposition}
\begin{proposition}\label{lemma-2}
For all $x\ge2$ satisfying (\ref{a0}) and (\ref{a2}), we have
\begin{equation}
\frac{\pr{(}\max_{1\leq k\leq n}\bar{S}_{k}^{N_0}\geq
x\bar{V}_{n}^{N_0} {)}}{1-\Phi(x)}
=2+C_{\delta,\tau}(\varepsilon^{-2}\Delta_{n,x}+\varepsilon).
\label{loser2}
\end{equation}
\end{proposition}

Indeed, noting that
\[
\frac x{\sqrt{2\uppi}(1+x^2)}\mathrm{e}^{-x^2/2}\le1-\Phi(x)
\le\frac1{\sqrt{2\uppi} x}\mathrm{e}^{-x^2/2}
\]
for $x\ge1$, we have that for $1\le k\le m=[x^2/2]$ and $x\ge1$,
\[
\frac{1-\Phi(\sqrt{x^2-k})}{1-\Phi(x)}\le2\mathrm{e}^{k/2}.
\]
This, together with (\ref{eq5})--(\ref{ad8}), implies that for all
$x\ge2$ satisfying (\ref{a0}) and (\ref{a2}),
\be\la{90}
&&\pr{\Bigl(\max_{1\leq k\leq n}S_{k}\geq xV_{n} \Bigr)} \no\\[-2pt]
&&\quad\leq \mathrm{e}^{-x^2}+2 \Biggl\{ 1-\Phi(x) +\sum_{k=1}^m\bigl\{1-\Phi\bigl(\sqrt
{x^2-k}\bigr)\bigr\}
\Biggl\{\sum_{j=1}^n\pr(|X_j|\ge\epsilon B_n/x) \Biggr\}^k \Biggr\} \no\\[-2pt]
&&\qquad\hphantom{\mathrm{e}^{-x^2}+}{} \times\{1+C_{\delta, \tau}
(\epsilon^{-2}\Delta_{n,x}+\epsilon) \} \\[-2pt]
&&\quad\leq 2\bigl(1-\Phi(x)\bigr) \{1+C_{\delta, \tau}
(\epsilon^{-3}\Delta_{n,x}+\epsilon+x^{-1}) \} \no\\[-2pt]
&&\quad\leq 2\bigl( 1-\Phi(x)\bigr) \bigl\{1+C_{\delta, \tau} \bigl(x^{-\min\{1/4,
\delta/20\}}+\Delta_{n,x}^{1/9} \bigr) \bigr\}.\no
\ee
Similarly, by
(\ref{ad6}), (\ref{ad8}) and (\ref{loser2}), we obtain that for
all $x\ge2$ satisfying (\ref{a0}) and~(\ref{a2}),
\be\la{91}
\pr{\Bigl(\max_{1\leq k\leq n}S_{k}\geq xV_{n} \Bigr)} \geq 2\bigl(1-\Phi(x)\bigr) \bigl\{1-C_{\delta, \tau}
(x^{-\min\{1/4, \delta/20\}}+\Delta_{n,x}^{1/9} ) \bigr\}.
\ee
Combining (\ref{90}) and (\ref{91}), we obtain (\ref{ad3}), and thus
Proposition~\ref{th2}.\vadjust{\goodbreak}

It remains to prove Propositions~\ref{lemma-1} and~\ref{lemma-2}, which
we give in
Sections~\ref{sec4} and~\ref{sec5}, respectively. The proof of
Theorem~\ref{th1} is now complete.
\end{pf*}
\section{\texorpdfstring{Proof of Proposition \protect\ref{lemma-1}}{Proof of Proposition 2}} \label{sec4}
Let $b=y/B^{N_l}_{n}$.
First, note that
\begin{eqnarray}\label{ad10}
\pr\Bigl(\max_{1\leq k\leq n}\bar{S}^{N_l}_{k}
\geq y\bar{V}^{N_l}_{n} \Bigr)
&\leq&\pr{\Bigl(2b\max_{1\leq k\leq
n}\bar{S}^{N_l}_{k}\geq(b\bar{V}^{N_l}_{n})^{2}+y^{2}
-\varepsilon^{2} \Bigr)}\nn
\\[-9pt]
\\[-9pt]
& &{} + \pr{\Bigl(\max_{1\leq k\leq
n}\bar{S}^{N_l}_{k}\geq y\bar{V}^{N_l}_{n}, |b\bar
{V}^{N_l}_{n}-y|\geq\varepsilon\Bigr)}.\nn
\end{eqnarray}
Furthermore, we have
\begin{eqnarray}\label{a15}
&&\pr{\Bigl(\max_{1\leq k\leq n}\bar{S}^{N_l}_{k}\geq y\bar
{V}_{n}^{N_l}, |b\bar{V}^{N_l}_{n}-y|\geq\varepsilon\Bigr)}
\nn\\
&&\quad \leq \pr{\Bigl(\max_{1\leq k\leq n}\bar{S}^{N_l}_{k}\geq
y\bar{V}_{n}^{N_l}, b^{2}(\bar{V}^{N_l}_{n})^{2}>y^{2}+\varepsilon y
\Bigr)}\nn
\\[-9pt]
\\[-9pt]
&&\qquad{} + \pr{\Bigl(\max_{1\leq k\leq n}\bar{S}^{N_l}_{k}\geq y\bar
{V}_{n}^{N_l}, b^{2}(\bar{V}^{N_l}_{n})^{2}<y^{2}-\varepsilon y \Bigr)}\nn\\
&&\quad =:  I_{1}+I_{2}\nn
\end{eqnarray}
and
\begin{eqnarray} \label{ad11}
&&\pr{\Bigl(2b\max_{1\leq k\leq n}\bar{S}^{N_l}_{k}\geq b^{2}
(\bar{V}_{n}^{N_l})^2+y^{2}-\varepsilon^{2} \Bigr)}\nn\\[-2pt]
&&\quad\leq
\pr{\Biggl(}\bigcup_{k=1}^{n} \{2b\bar{S}^{N_l}_{k}
\geq b^{2}(\bar{V}_{n}^{N_l})^2+ y^{2}-\varepsilon^{2},\nn\\[-2pt]
&&
\qquad\hphantom{\pr{\Biggl(}\bigcup_{k=1}^{n} \{}\ep[(\bar{V}_{n}^{N_l})^2-(\bar{V}_{k}^{N_l})^2 ]
- [(\bar{V}_{n}^{N_l})^2-(\bar{V}_{k}^{N_l})^2 ]
\geq
\varepsilon^{2}/b^2 \} {\Biggr)}\\[-2pt]
&&\qquad{}+\pr{\Biggl(\bigcup_{k=1}^{n} \{2b\bar{S}^{N_l}_{k}\geq
b^{2}(\bar{V}_{k}^{N_l})^2+b^{2}\ep[
(\bar{V}_{n}^{N_l})^2-
(\bar{V}_{k}^{N_l})^2 ]+y^{2}-2\varepsilon^{2} \} \Biggr)}\nn\\[-2pt]
&&\quad =:I_{3}+I_{4}.\nn
\end{eqnarray}
By (\ref{ad10})--(\ref{ad11}),
Proposition~\ref{lemma-1} follows from the following Lemma~\ref{1}.
\begin{lemma} \la{1} Under the conditions of Proposition \ref
{lemma-1}, we have
\be
\label{ad12} I_1 &\le& C_{\delta, \tau} y^{-2}\exp(-y^{2}/2), \\[-2pt]
\label{ad12-1}I_2&\le& C_{\delta, \tau} y^{-2}\exp(-y^{2}/2),  \\[-2pt]
\label{ad12-2}I_3&\le& C_{\delta, \tau} y^{-2}\exp(-y^{2}/2),  \\[-2pt]
\label{ad14}I_4 &\le& 2[1-\Phi(y)]
[1+C_{\delta,\tau}(\varepsilon^{-2}\Delta_{n,x}+\varepsilon) ].
\ee
\end{lemma}

To prove Lemma~\ref{1}, we start with some preliminaries.
Note that
\be\gamma\max\{ x^{-1/2}, x^{-\delta/10}\}\le
\varepsilon\le\min\{1/24,\delta/72\},\qquad \Delta_{n,x}\le
(\varepsilon/2)^{9/2}.\label{bor} \ee
This fact (\ref{bor}) is repeatedly used in the proof without further
explanation. Define
$k_0=0, k_T=n$ and $k_{i}, 1\le i<T$, by
\[
k_{i}=\max\Biggl\{k\dvt \sum_{j=k_{i-1}+1}^{k}\ep X^{2}_{j}\leq2^{-1}
\varepsilon^{3} B^{2}_{n}/x^{2} \Biggr\}.
\]
By the definition of $k_{i}$,
\begin{equation}\lbl{new-1}
\sum_{j=k_{i-1}+1}^{k_{i}}\ep X^{2}_{j}\leq2^{-1}\varepsilon^{3}
B^{2}_{n}/x^{2} \quad\mbox{and}\quad \sum_{j=k_{i-1}+1}^{k_{i}+1}\ep X^{2}_{j}>
2^{-1}\varepsilon^{3} B^{2}_{n}/x^{2}
\end{equation}
for any $1\le i<T$. By (\ref{a0}) and (\ref{bor}),
\begin{eqnarray}\label{a1}
x^{2}\max_{1\leq k\leq n}\ep X^{2}_{k}&\leq& x^2 \max_{1\le k\le n} [\ep
X^{2}_{k}\{|X_{k}|>B_{n}/x\}+(\ep|X_{k}|^{3}I\{|X_{k}|\leq B_{n}/x\}
)^{2/3} ] \nonumber
\\[-8pt]
\\[-8pt]
&\le& B_n^2(\Delta_{n,x}+\Delta_{n,x}^{2/3}) \le\varepsilon^{3}
B^{2}_{n}/4 ,\no
\end{eqnarray}
which, together with \eq{new-1}, implies that
\[
\sum_{j=k_{i-1}+1}^{k_{i}}\ep X^{2}_{j}\geq4^{-1}\varepsilon^{3}
B^{2}_{n}/x^{2}.
\]
Therefore,
\[
(T-1)4^{-1}\varepsilon^{3} B^{2}_{n}/x^{2}\leq
\sum_{i=1}^{T-1}\sum_{j=k_{i-1}+1}^{k_{i}}\ep X^{2}_{j}\leq B^{2}_{n},
\]
which yields $T\leq4x^{2}/\varepsilon^{3}+1$. For
$ k_{i-1}+1\leq j\leq k_{i}-1$, define events
\[
\A_{j}=\bigl \{\bar{S}^{N_l}_{j}\geq y\sqrt{(B^{N_l}_{n})^2(1+
\varepsilon/y)} \bigr\},\qquad \C_{j}= \Biggl\{\sum_{k=j+1,k\notin N_l}^{k_{i}}(\bar
{X}_{k}-\ep
\bar{X}_{k})\geq-\varepsilon B^{N_l}_{n}/y \Biggr\}.
\]
Note that $\sum_{k\in
N_l}\ep X^{2}_{k}\le\varepsilon^3B_n^2/8$ for all $0\le
l\le m=[x^2/2]$ by (\ref{a1}),
and thus
\be
B_n^2\ge(B^{N_l}_{n})^{2}=
\sum_{k=1}^{n}\ep X^{2}_{k}-\sum_{k\in N_l}\ep
X^{2}_{k}\geq(1-\varepsilon^{3}/8)B^{2}_{n}\ge\frac78B_n^2. \label{ade}
\ee
Applying the Chebyshev inequality, we have, for any $k_{i-1}\leq j\leq
k_{i}$ and $x/2\leq y\le x$,
\begin{equation}\label{har1}
\pr(\C_{j})\geq1-\frac{y^{2}\sum_{k=j+1}^{k_{i}}\ep X^{2}_{k}}
{\varepsilon^{2}(B^{N_l}_{n})^{2}}\geq
1-4\varepsilon/7\ge1/2 .
\end{equation}

We are now ready to prove Lemma~\ref{1}.
\begin{pf*}{Proof of (\ref{ad12})} It follows from
(\ref{har1}) and the independence between $\C_j$ and $\{\A_l, l \leq
j\}$ that
\begin{eqnarray}\label{ad19}
I_{1}&\leq& \sum_{i=1}^T \pr\Biggl(\bigcup_{j=k_{i-1}+1}^{k_i}\A_{j} \Biggr) \no\\[-2pt]
&\leq&\sum_{i=1}^{T} {\Biggl[\pr(\A_{k_{i-1}+1})+\sum
_{j=k_{i-1}+2}^{k_{i}}\pr{(}\A^{c}_{k_{i-1}},\ldots, \A
^{c}_{j-1},\A_{j} {)} \Biggr]}\nn
\\[-9pt]
\\[-9pt]
&\leq& 2 \sum_{i=1}^{T} {\Biggl[\pr(\A_{k_{i-1}+1},
\C_{k_{i-1}+1})+\sum_{j=k_{i-1}+2}^{k_{i}}\pr{(}\A
^{c}_{k_{i-1}},\ldots,
\A^{c}_{j-1}, \A_{j}, \C_{j} {)} \Biggr]}\nn\\[-2pt]
&\leq&
2 \sum_{i=1}^{T}\pr{\bigl(
\bar{S}^{N_l}_{k_{i}}-\ep\bar
{S}^{N_l}_{k_{i}}\geq
y\sqrt{(B^{N_l}_{n})^2(1+\varepsilon
/y)}-\varepsilon B^{N_l}_{n}/y-D_{k_{i}} \bigr)},\nn
\end{eqnarray}
where $D_{k_{i}}=\sum_{j=1}^{k_{i}}\ep|X_{j}|I\{|X_{j}|> \varepsilon
B_{n}/x\}$.
Taking $t=y\sqrt{1+\varepsilon/y}/B^{N_l}_{n}$
and noting
\bestar
t(\varepsilon B_{n}/y+D_{k_{i}})\leq2\varepsilon+1,
\eestar
we have
\begin{eqnarray*}
&&\pr{\bigl(\bar{S}^{N_l}_{k_{i}}-\ep
\bar{S}^{N_l}_{k_{i}}\geq
y\sqrt{(B^{N_l}_{n})^2(1+\varepsilon/y)}-\varepsilon
B_n^{N_l}/y- D_{k_{i}}\bigr )}\\[-2pt]
&& \quad\leq
9\exp(-y^{2}-\varepsilon y)\prod_{j=1,j\notin
N_l}^{k_{i}}\ep\exp\bigl(t(\bar{X}_{j}-\ep\bar{X}_{j}) \bigr)\\[-2pt]
&&\quad \leq9\exp(-y^{2}-\varepsilon y) \prod_{j=1,j\notin
N_l}^{k_{i}} {\biggl(1+\frac{\ep X^{2}_{j}}{2}t^{2}+8
t^{3} \ep|\bar{X}_{j}|^{3}\mathrm{e}^{2t\varepsilon B_n/x} \biggr)}\\[-2pt]
&&\quad\leq9 \exp(-y^{2}/2-\varepsilon y/2+A\Delta_{n,x}).
\end{eqnarray*}
Submitting this estimate into (\ref{ad19}) and
recalling $T\le4x^2/\varepsilon^3+1$, $x/2\leq y\le x$ and
$\varepsilon\ge\gamma x^{-1/2}$,
we obtain
\begin{eqnarray}\label{har2}
I_{1}&\leq& (4\varepsilon^{-3}x^{2}+1)\exp(-y^{2}/2-\varepsilon y/2
+A\Delta_{n,x}) \no
\\[-9pt]
\\[-9pt]
&\le&C_{\delta, \tau} y^{-2} \mathrm{e}^{-y^2/2}. \no
\end{eqnarray}
This proves (\ref{ad12}).
\end{pf*}
\begin{pf*}{Proof of (\ref{ad12-1})} For this part,
let $Y_{k_{i}}=\sum_{j=k_{i-1}+1,j\notin
N_l}^{k_{i}}\bar{X}^{2}_{j}$, and define
\[
\bar{\A}_{j}=\bigl\{\bar{S}^{N_l}_{j}\geq y
\sqrt{(\bar{V}_{n}^{N_l})^2-Y_{k_{i}}}, b^{2} [(\bar
{V}_{n}^{N_l})^2-Y_{k_{i}} ]<y^{2}-
\varepsilon y\bigr\}, \qquad1\leq j\leq n.\vadjust{\goodbreak}
\]
From (\ref{har1}) and
the independence between $\C_j$ and $\{\bar{\A}_{l}, l \leq j\}$, it
follows that
\begin{eqnarray}\label{har4}
I_{2}&\leq& \sum_{i=1}^{T} {\Biggl[\pr(\bar{\A}_{k_{i-1}+1})+
\sum_{j=k_{i-1}+2}^{k_{i}}\pr{(}\bar{\A}^{c}_{k_{i-1}},\ldots,
\bar{\A
}^{c}_{j-1},
\bar{\A}_{j} {)} \Biggr]}\nn\\[-2pt]
&\leq&2\sum_{i=1}^{T} {\Biggl[\pr(\bar{\A}_{k_{i-1}+1},\C
_{k_{i-1}+1})+\sum
_{j=k_{i-1}+2}^{k_{i}}\pr{(}\bar{\A}^{c}_{k_{i-1}},\ldots,
\bar{\A}^{c}_{j-1}, \bar{\A}_{j},\C_{j} {)} \Biggr]}\nn
\\[-9pt]
\\[-9pt]
&\leq&2\sum
_{i=1}^{T}\pr{\bigl(\bar{S}^{N_l}_{k_{i}}-\ep
\bar{S}^{N_l}_{k_{i}}\geq y\sqrt{(\bar
{V}_{n}^{N_l})^2-Y_{k_{i}}}-\varepsilon B^{N_l}_{n}/y-D_{k_{i}},
b^{2} [(\bar{V}_{n}^{N_l})^2-Y_{k_{i}} ]<y^{2}-\varepsilon y
\bigr)}\nn\\[-2pt]
&=:& 2\sum_{i=1}^{T}I_{2,i}, \nn
\end{eqnarray}
where, as before, $D_{k_{i}}=\sum_{j=1}^{k_{i}}\ep|X_{j}|I\{|X_{j}|>
\varepsilon B_{n}/x\}$. Furthermore, for $i=1, \ldots,T$,
\begin{eqnarray*}
I_{2,i} &\leq&
\pr{\bigl((\bar{V}_{n}^{N_l})^2-Y_{k_{i}}<(1-\varepsilon)
(B^{N_l}_{n})^2 \bigr)}\\[-2pt]
& &{}+\sum_{k=1}^{[y]}
\pr{\bigl(}\bar{S}^{N_l}_{k_{i}}-\ep
\bar{S}^{N_l}_{k_{i}}\geq y
\sqrt{(B^{N_l}_{n})^2 [1-(k+1)\varepsilon/y ]}-\varepsilon
B^{N_l}_{n}/y-D_{k_{i}},\\[-2pt]
& &\hphantom{{}+\sum_{k=1}^{[y]}
\pr{\bigl(}}
(B^{N_l}_{n})^2 [1-(k+1)\varepsilon/y ]
<(\bar{V}_{n}^{N_l})^2-
Y_{k_{i}}<(B^{N_l}_{n})^2 [1-k\varepsilon
/y ] {\bigr)} \\[-2pt]
&=:&I_{2,i,0}+ \sum_{k=1}^{[y]}I_{2,i,k}.
\end{eqnarray*}
Note that, for any $t_1\ge0$ and $t_2\ge0$,
\begin{eqnarray} \label{har5}
&&\ep\exp\bigl(t_{1}(\bar{X}_{k}-\ep\bar{X}_{k})+t_{2}(\ep\bar{X}^{2}_{k}-
\bar{X}^{2}_{k}) \bigr) \no\\
&&\quad \leq 1+\tfrac12 \ep\bigl( t_{1}(\bar{X}_{k}-\ep\bar
{X}_{k})+t_{2}(\ep\bar{X}^{2}_{k}-
\bar{X}^{2}_{k}) \bigr)^2\no\\
& &\qquad{} + ( 8t^{3}_{1}\ep|\bar{X}_{k}|^{3} +
8t_2^3 \ep|\bar{X}_{k}|^{6} )
\mathrm{e}^{2t_{1}\varepsilon
B_{n}/x+t_{2}\ep X^{2}_{k}}
\\
& &\quad\leq \exp{\bigl(}\tfrac12 t^{2}_{1} \ep\bar{X}^{2}_{k}+ \tfrac12
(4t_{1}t_{2}+t^{2}_{2}\varepsilon
B_{n}/x )\ep|\bar{X}_{k}|^{3}\no\\
&&\qquad\hphantom{\exp{\bigl(}}{} +(8t^{3}_{1}+8t_2^3\varepsilon^{3}B_{n}^3/x^{3})\ep|\bar{X}_{k}|^{3}
\mathrm{e}^{2t_{1}\varepsilon
B_{n}/x+t_{2}\max_{1\leq k\leq n}\ep X^{2}_{k}} {\bigr)}.\no
\end{eqnarray}
Let $t_{1}=y\sqrt{1-(k+1)\varepsilon/y}/B_{n}^{N_l}$ and
$t_{2}=\varepsilon^{-1}y^{2}/(B_{n}^{N_l})^2$ in (\ref{har5}). Noting that
\[
t_1(\varepsilon B_n^{N_l}/y+D_{k_i}) \le\varepsilon+1
,
\]
we have for $1\le k\le[x]$,
\begin{eqnarray*}
I_{2,i,k} &\leq& \pr{\bigl(}t_{1}(\bar{S}^{N_l}_{k_{i}} -\ep
\bar{S}^{N_l}_{k_{i}})+t_{2} \{\ep
[(\bar{V}_{n}^{N_l})^2-Y_{k_{i}}]-[(\bar{V}_{n}^{N_l})^2
-Y_{k_{i}}] \}\\
&& \quad\geq
y^{2}-(k+1)\varepsilon/y+t_{2}k\varepsilon
(B^{N_l}_{n})^2/y -2 {\bigr)}\\
& \leq&
\exp{\bigl(-y^{2}+ (k+1)\varepsilon/y-t_{2}k\varepsilon
(B^{N_l}_{n})^2/y+ 2 \bigr)}\\
&& {}\times
\prod_{k=1,k\notin
N_l}^{k_{i-1}}\ep\exp{\bigl(t_{1}(\bar{X}_{k}-\ep
\bar{X}_{k}) +t_{2}(\ep\bar{X}^{2}_{k}-\bar{X}^{2}_{k} )\bigr )}\\
&&{} \times\prod_{k=k_{i-1}+1,k\notin
N_l}^{k_{i}}\ep\exp{\bigl(t_{1}(\bar{X}_{k}-\ep
\bar{X}_{k})\bigr )} \times
\prod_{k=k_{i}+1,k\notin
N_l}^{n}\ep\exp{\bigl(t_{2}(\ep\bar{X}^{2}_{k}-\bar{X}^{2}_{k}
) \bigr)}\\
& \leq& \exp\Biggl(-y^{2}/2+
2^{-1}(k+1)\varepsilon/y-t_{2}k\varepsilon
(B^{N_l}_{n})^2/y +2\cr
&&\hphantom{\exp{\Biggl(}}{}+A {(}t_{1}t_{2}+ t^{2}_{2}\varepsilon
B_{n}/x +t^{3}_{1}+t^{3}_{2}\varepsilon^{3} B^{3}_{n}/x^{3}
{)}\sum_{k=1}^{n}\ep|\bar{X}_{k}|^{3} \Biggr)\\
&\leq&
\exp\bigl(-y^{2}/2+ 2^{-1}(k+1)\varepsilon/y-ky+
A\varepsilon^{-1}\Delta_{n,x}+2 \bigr)\\
& \leq&
A \exp(-y^{2}/2-y/2).
\end{eqnarray*}
Similarly, by (\ref{har5}) with $t_1=0$, we have
\begin{eqnarray*}
I_{2,i,0} &\leq& \pr{\bigl(}t_{2} \{\ep[(\bar{V}_{n}^{N_l})^2
-Y_{k_{i}} ]- [(\bar{V}_{n}^{N_l})^2-Y_{k_{i}} ] \}\geq
t_{2}\varepsilon
(B^{N_l}_{n})^2-\varepsilon^{2} {\bigr)}\\
&\leq& A\exp(-y^{2}+A\varepsilon^{-1}\Delta_{n,x})\le
A_1\exp(-y^2/2-y).
\end{eqnarray*}
Combining above inequalities yields
\begin{equation}\label{a17}
I_{2}\leq A(4x^2/\varepsilon^3+1) \mathrm{e}^{-y^{2}/2-y}\le A_1
y^{-2}\mathrm{e}^{-y^2/2}.
\end{equation}
The proof of (\ref{ad12-1}) is now complete.
\end{pf*}
\begin{pf*}{Proof of (\ref{ad12-2})}
Following the arguments in the estimates of $I_1$ and $I_2$, we have
\begin{eqnarray}
I_3&\leq&\sum_{i=1}^{T}
\pr{\Biggl(}\bigcup_{j=k_{i-1}+1}^{k_{i}}\{2b\bar{S}^{N_l}_{j}\geq
b^{2}(\bar{V}_{n}^{N_l})^2+ y^{2}-\varepsilon^{2},\no\\
&&\hphantom{\sum_{i=1}^{T}
\pr{\Biggl(}\bigcup_{j=k_{i-1}+1}^{k_{i}}\{}
\ep[(\bar{V}_{n}^{N_l})^2-(\bar{V}^{N_l}_{j})^2 ]-
[(\bar{V}_{n}^{N_l})^2-(\bar{V}^{N_l}_{j})^2
]\geq
\varepsilon^{2}(B^{N_l}_{n})^2/y^{2}\} {\Biggr)}\hspace*{30pt}\\
&\leq&
\sum_{i=1}^{T}
\pr{\Biggl(}\bigcup_{j=k_{i-1}+1}^{k_{i}}\{2b\bar{S}^{N_l}_{j}\geq
b^{2} [(\bar{V}_{n}^{N_l})^2-Y_{k_{i}} ] +
y^{2}-\varepsilon^{2}\},\no \\
&& \hphantom{\sum_{i=1}^{T}
\pr{\Biggl(}}\sum_{k=k_{i}+1,k\notin
N_l}^{n}(\ep\bar{X}^{2}_{k}-\bar{X}^{2}_{k})\geq
2^{-1}\varepsilon^{2}(B^{N_l}_{n})^2/y^{2} {\Biggr)}\no\\
&\leq& 2\sum_{i=1}^{T}
\pr{\Biggl(}2b\bar{S}^{N_l}_{k_{i}}\geq
b^{2} [(\bar{V}_{n}^{N_l})^2-Y_{k_{i}} ]+
y^{2}-2\varepsilon, \no\\
&& \hphantom{2\sum_{i=1}^{T}
\pr{\Biggl(}}\sum_{k=k_{i}+1,k\notin
N_l}^{n}(\ep\bar{X}^{2}_{k}-\bar{X}^{2}_{k})\geq
2^{-1}\varepsilon^{2}(B^{N_l}_{n})^2/y^{2} {\Biggr)} \no\\
&=:& 2\sum_{i=1}^TI_{3i}.\no
\end{eqnarray}
As in the proof of (\ref{har5}), it can be easily shown that for
$\al\geq0$,
\be\label{wed1}
\ep \mathrm{e}^{b\bar X_j-\al b^2\bar X_j^2} \le
\exp\bigl\{(1/2-\al)b^2\ep\bar X_j^2+A\Delta_{n,x}^{(j)}\bigr\},
\ee
where
\bestar
\Delta_{n,x}^{(j)} &=& \frac{x^2}{B_n^2}\ep
X_j^2I\{|X_j|\ge\varepsilon B_n/x\}+
\frac{x^3}{B_n^3}\ep|X_j|^3I\{|X_j|\le\varepsilon B_n/x\} \no\\
&\le& \varepsilon^{-1} \biggl(\frac{x^2}{B_n^2}\ep X_j^2I\{|X_j|\ge
B_n/x\}+ \frac{x^3}{B_n^3}\ep|X_j|^3I\{|X_j|\le B_n/x\} \biggr)
\eestar
and
\begin{eqnarray}\label{wed2}
&&\ep \mathrm{e}^{\al(\ep\bar X_j^2-\bar X_j^2)-b^2\bar X_j^2/2} \nn
\\[-8pt]
\\[-8pt]
&&\quad\le \exp\bigl\{-\tfrac12b^2\ep\bar X_j^2 +(2\al^2B_n/x+x^3/B_n^3)
\varepsilon\ep|\bar X_j|^3\mathrm{e}^{\al
\max_{1\leq k\leq n}\ep X^{2}_{k}}\bigr \}. \nn
\end{eqnarray}
Next, let
$t$ satisfy
\[
t\mathrm{e}^{t\max_{1\leq k\leq n}\ep X^{2}_{k}} =\frac{\varepsilon B_{n}}{24x
\sum_{j=k_{i}+1}^{n}\ep|\bar{X}_{j}|^{3}}.
\]
Clearly $t$ exists.
Furthermore, we have $t\ge x^2/B_n^2$. Indeed,
if $t\max_{1\leq k\leq n}\ep X^{2}_{k}\geq\varepsilon$, then by
(\ref
{a1}) and recalling $\varepsilon\le1/24$,
\[
t\ge\varepsilon/\max_{1\leq k\leq n}\ep X^{2}_{k}\ge
4\varepsilon^{-2}x^2/B_n^2 \ge x^2/B_n^2.
\]
If
$t\max_{1\leq k\leq n}\ep X^{2}_{k}\leq\varepsilon$, then
\bestar
t\ge \frac{\varepsilon B_{n}}{24\mathrm{e}^{\varepsilon}
x\sum_{j=k_{i}+1}^{n}\ep|\bar{X}_{j}|^{3}} \ge\frac{\varepsilon
x^2}{30 B_n^2 \Delta_{n,x}} \ge\frac1{15}
\Delta_{n,x}^{-7/9} x^2/B_n^2 \ge x^2/B_n^2.
\eestar
Now it follows
from (\ref{wed1}) and (\ref{wed2}) with $\al=t$ that
\begin{eqnarray}\label{an1}
 I_{3i}&\leq&
\pr{\Biggl(}b\bar{S}^{N_l}_{k_{i}}- 2^{-1}b^{2} [(\bar
{V}_{n}^{N_l})^2-Y_{k_{i}} ]+
t\sum_{k=k_{i}+1,k\notin
N_l}^{n}(\ep\bar{X}^{2}_{k}-\bar{X}^{2}_{k})\no\\
&&\hphantom{\pr{\Biggl(}} \geq
y^{2}/2-2\varepsilon+
2^{-1}t\varepsilon^{2}(B^{N_l}_{n})^2/y^{2} {\Biggr)} \no\\
&\le&
\exp\biggl[2\varepsilon-y^{2}/2-\frac{\varepsilon^{2}(B^{N_l}_{n})^2t}{2y^{2}} \biggr] \prod_{j=1,j\notin
N_l}^{k_{i-1}}\ep
\mathrm{e}^{b\bar{X}_{j}-2^{-1}b^{2}\bar{X}^{2}_{j}} \no\\
&&{} \times\prod_{j=k_{i-1}+1,j\notin N_l}^{k_{i}}\ep \mathrm{e}^{b\bar
{X}_{j}} \times\prod_{j=k_{i}+1,j\notin N_l}^{n}\ep
\mathrm{e}^{-2^{-1}b^{2}\bar{X}^{2}_{j}+t(\ep
\bar{X}^{2}_{j}-\bar{X}^{2}_{j})}\no
\\[-8pt]
\\[-8pt]
&\leq&
A \exp(-y^{2}/2)\exp{\Biggl(}\varepsilon^{-1}\Delta_{n,x}-
\frac{\varepsilon^{2}B^{2}_{n}t}{3x^{2}}-\frac{y^{2}\sum
_{j=k_{i}+1}^{n}\ep
X^{2}_{j}}{2(B_{n}^{N_l})^2}+\frac{y^{2}\sum_{j=k_{i-1}+1}^{k_i}\ep
X^{2}_{j}}{2(B_{n}^{N_l})^2}\hspace*{32pt}\no \\
&&\hphantom{A \exp(-y^{2}/2)\exp{\Biggl(}}{} +
(2t^2B_n/x+x^3/B_n^3) \varepsilon\sum_{j=k_{i}+1}^{n}\ep|\bar
X_j|^3\mathrm{e}^{t
\max_{1\leq k\leq n}\ep X^{2}_{k}} {\Biggr)}\no\\
&\leq&
A_1 \exp(-y^{2}/2)\exp\biggl(-\frac{\varepsilon^{2}B^{2}_{n}t}{4x^{2}}-
\frac{y^{2}\sum_{j=k_{i}+1}^{n}\ep X^{2}_{j}}{4B_{n}^2} +\frac{
\varepsilon^2x^2}{12B_n^2 t}
\biggr) \no\\
&\leq&
A_1 \exp(-y^{2}/2)\exp\biggl(-\frac{\varepsilon^{2}B^{2}_{n}t}{4x^{2}}-
\frac{y^{2}\sum_{j=k_{i}+1}^{n}\ep X^{2}_{j}}{4B_{n}^2} \biggr) .\no
\end{eqnarray}
Note that when $t\le\frac{2\delta x^2\log x}{B_n^2\varepsilon^3}$,
$t\max_{1\leq k\leq}\ep X^{2}_{k}< \frac{\delta}{2}\log x$ by (\ref
{a1}). Thus, by the
definition of $t$,
\[
t\geq\frac{\varepsilon B_{n}}{24 x^{1+\delta/2}
\sum_{j=k_{i}+1}^{n}\ep|\bar{X}_{j}|^{3}}.
\]
Now considering $t\le\frac{2\delta x^2\log x}{B_n^2\varepsilon^3}$
and $t \geq\frac{2\delta x^2\log x}{B_n^2\varepsilon^3}$,
we have, by \eq{an1},
\be\label{hope}
I_{3i} &\le& A y^{-\delta/(3\varepsilon)}\exp(-y^{2}/2) \no
\\[-8pt]
\\[-8pt]
&&{} + A \mathrm{e}^{-y^{2}/2} \exp\biggl{(}-\frac{y^2\sum_{j=k_{i}+1}^{n}\ep
X^{2}_{j}}{4B^{2}_{n}}-
\frac{\varepsilon^{3}B^{3}_{n}}{144x^{3+\delta/2}
\sum_{j=k_{i}+1}^{n}\ep|\bar{X}_{j}|^{3}} \biggr{)}.\no
\ee
From the definition of $n_0$,
$\sum_{j=n_{0}+1}^{n}\ep X^{2}_{j}\leq192B^{2}_{n}x^{-2}\log x$ and
thus by (\ref{a2})
\[
\sum_{j=n_{0}+1}^{n}\ep|\bar{X}_{j}|^{3}\leq\frac{192\tau
B^{3}_{n}\log x}{x^{3+\delta}}.
\]
For $i<i_0$, where $i_{0}=\max\{i\dvt k_{i}+1\leq n_{0}\}$, we have
\[
y^2\sum_{j=k_i}^{n}\ep X^{2}_{j} \geq x^2\sum_{j=n_0}^{n}\ep
X^{2}_{j}/4\ge
24 B_n^2\log x.
\]
It now follows from \eq{hope}, \eq{bor} and the fact $T\le
4x^2/\varepsilon^3+1$ that
\begin{eqnarray}\label{a14}
I_3 &\le& 2\sum_{i=1}^TI_{3i} \no\\
&\le& 2A T y^{-\delta/(3\varepsilon)}\mathrm{e}^{-y^2/2}+ 2A \mathrm{e}^{-y^2/2} i_0
\mathrm{e}^{-6\log x}\no\\
&&{} +2A \mathrm{e}^{-y^2/2}
\sum_{i=i_{0}+1}^{T}\exp\biggl{(}-\frac{\varepsilon
^{3}B^{3}_{n}}{144x^{3+\delta/2}\sum_{j=k_i+1}^{n}\ep|\bar{X}_{j}|^{3}}
\biggr{)}\\
&\leq& A_1 (4x^2/\varepsilon^3+1)\mathrm{e}^{-y^2/2} (y^{-6}+
\mathrm{e}^{-A x^{\delta/2}\varepsilon^3/\log x} )\no\\
&\leq&
C_{\delta, \tau}y^{-2}\mathrm{e}^{-y^2/2}.\no
\end{eqnarray}
This completes the proof of \eq{ad12-2}.
\end{pf*}
\begin{pf*}{Proof of (\ref{ad14})} For this result, we need the following
moderate deviation theorem for the standardized sum
due to Sakhanenko \cite{Sak92} (also see Heinrich \cite{H}).
\begin{lemma}\label{le3}
Suppose that $\eta_{1},\ldots,\eta_{n}$ are independent random
variables such that $\ep\eta_{j}=0$ and $|\eta_j|\le1$ for $j\geq
1$. Write $\sigma_n^{2}=\sum_{j=1}^{n}\ep\eta^{2}_{j}$ and
${\mathcal
L}_n=\sum_{j=1}^{n}\ep|\eta_{j}|^{3}/\sigma_n^3$. Then there exists
an absolute constant $A>0$ such that for all $1\le x\le\min\{\sigma_n,
{\mathcal L}_n^{-1/3}\}/A$,
\begin{equation}\label{dis}
\frac{\pr(\sum_{j=1}^{n}\eta_{j}\geq
x\sigma_n)}{1-\Phi(x)}=1+\mathrm{O}(1)x^3{\mathcal L}_n,
\end{equation}
where $|\mathrm{O}(1)|$ is bounded by an absolute constant.
\end{lemma}

To prove (\ref{ad14}), write
\begin{eqnarray*}
\xi_{j}&=&2b\bar{X}_{j}-b^{2}\bar{X}^{2}_{j}+b^{2}\ep\bar
{X}^{2}_{j},\cr
\E_{j}&=&\Biggl\{\sum_{k=1,k\notin N_l}^{j}\xi_{k}\geq
z\Biggr\}, \qquad\mbox{where } z=2(y^{2}-\varepsilon^{2}).
\end{eqnarray*}
Note that $|\xi_j-\ep\xi_j|\le4\varepsilon+2\varepsilon^2\le
5\varepsilon$, and by the non-uniform Berry--Esseen bound, there
exists an absolute constant $A_0$ such that for any $1\leq k\leq n$
and $c>0$,
\begin{eqnarray*}
&&\pr\Biggl{(}\sum_{j=k,j\notin N_l}^{n}(\xi_{j}-\ep\xi_{j}) \leq
-c\varepsilon\Biggr{)} \\
&&\quad\leq
1-\Phi(t)+\frac{A_{0}\sum_{j=k,j\notin N_l}^{n}\ep|\xi
_{j}-\ep\xi_j|^{3}}{(1+t)^3s^{3}_{n,k}}\\
&&\quad\leq\frac{1}{2}-\frac{1}{\sqrt{2\uppi}}\int
_{0}^{t}\mathrm{e}^{-s^{2}/2}\,\mathrm{d}s+\frac
{5A_{0}}{c}(1+t)^{-3}t,
\end{eqnarray*}
where $s^{2}_{n,k}=\sum_{j=k,j\notin N_l}^{n}\Var(\xi_{j})$
and $t=c\varepsilon/s_{n,k}$.
Because
$\int_{0}^{t}\mathrm{e}^{-s^{2}/2}\,\mathrm{d}s\geq t(1+t)^{-3}/2$ for any $t\geq0$, we
may choose $c_{0}\geq10 A_{0}\sqrt{2\uppi}$ such that
for all $1\le k\le n$,
\begin{equation}\label{a18}
\pr\Biggl{(}\sum_{j=k,j\notin N_l}^{n}(\xi_{j}-\ep\xi_{j}) \leq
-c_0\varepsilon\Biggr{)}\leq 1/2.
\end{equation}
By virtue of (\ref{a18}), we obtain that
\begin{eqnarray}\label{a19}
I_{4}&=&\pr(\E_{1})+\sum_{k=2}^{n}\pr(\E^{c}_{1},\ldots, \E
^{c}_{k-1},\E
_{k})\no\\
&\leq&2\pr\Biggl(\E_{1}, \sum_{j=2,j\notin N_l}^{n}(\xi_{j}-\ep\xi
_{j}) \geq-c_{0}\varepsilon\Biggr)\no
\\[-8pt]
\\[-8pt]
& &{}+2\sum_{k=2}^{n}\pr\Biggl(\E^{c}_{1},\ldots, \E^{c}_{k-1},\E
_{k},\sum
_{j=k+1,j\notin N_l}^{n}(\xi_{j}-\ep\xi_{j}) \geq
-c_{0}\varepsilon\Biggr)\no\\
&\leq&2\pr\Biggl{(}\sum_{k=1,k\notin N_l}^{n}(\xi_{k}-\ep\xi
_{k})\geq z-c_{0}\varepsilon-D_{n} \Biggr{)},\no
\end{eqnarray}
where $D_{n}=\sum_{j=1, j\notin N_l}^{n}|\ep\xi_j|$. Write
$z'=z-c_{0}\varepsilon-D_{n}$.
It is not difficult to show that
\bestar
D_n &\le& 2b \sum_{j=1, j\notin N_l}^n\ep|X_j|I\{|X_j|\ge
\varepsilon B_n/x\}\le4\varepsilon^{-2}\Delta_{n,x},\\
s^{2}_{n,1} &=&\sum_{j=1,j\notin N_l}^{n}\Var(\xi_{j}) = 4b^2
\sum_{j=1, j\notin N_l}^n\ep X_j^2+\mathrm{O}(1) \varepsilon^{-1}
\Delta_{n,x} \\
&=& 4y^2+\mathrm{O}(1) \varepsilon^{-2} \Delta_{n,x},
\eestar
where
$|\mathrm{O}(1)|\le30$. This yields that
\bestar
\frac{z'}{s_{n,1}} &= &
y+\mathrm{O}(1) [(\varepsilon+\varepsilon^{-2}\Delta_{n,x})/y ],
\eestar
 where $|\mathrm{O}(1)|\le40$. Therefore, by Lemma \ref
{le3} with $\eta_j=\xi_j-\ep\xi_j$
\be\label{a21}
I_4 &\le& 2 [1-\Phi(z'/s_{n,1}) ] \Biggl[1+A (z'/s_{n,1})^3
s_{n,1}^{-3}\sum
_{j=1, j\notin N_l}^n \ep|\xi_j|^3 \Biggr] \no
\\[-8pt]
\\[-8pt]
&\le& 2
[1-\Phi(y) ] [1+A (\varepsilon+\varepsilon^{-2}\Delta_{n,x}) ],\no
\ee
where we have used the fact that whenever $x
\theta_n\to0$, \bestar\frac
{1-\Phi(x+\theta_n)}{1-\Phi(x)}=1+\mathrm{O}(1)x \theta_n. \eestar This
proves (\ref{ad14}), and also completes the proof of Proposition
\ref{lemma-1}.
\end{pf*}
\section{\texorpdfstring{Proof of Proposition \protect\ref{lemma-2}}{Proof of Proposition 3}} \label{sec5}
By Proposition~\ref{lemma-1}, it suffices to show that
\be\label{los9}
\pr\Bigl{(}\max_{1\leq k\leq n}\bar{S}_{k}\geq x\bar{V}_{n} \Bigr{)}
\geq2 \bigl(1-\Phi(x) \bigr) \bigl(1-C_{\delta, \tau}
(\varepsilon^{-2}\Delta_{n,x}+\varepsilon)\bigr ).
\ee
Toward this end, let $b=x/B_{n}^{N_0}$ throughout this
section. Recall
(\ref{bor}), which we use repeatedly in the proof without further
explanation. Let
$n_0$ be defined as in (\ref{n0a}).
It can be readily seen that
\begin{eqnarray}\label{bad1}
&&\pr\Bigl{(}\max_{1\leq k\leq n}\bar{S}_{k}\geq x\bar{V}_{n} \Bigr{)}\nn\\[-2pt]
&&\quad \geq \pr\Bigl{(}2b\max_{n_0\leq k\leq n}\bar{S}_{k}\geq
b^{2}\bar{V}^{2}_{n}+x^{2} \Bigr{)}\nn\\[-2pt]
&&\quad\geq \pr\Biggl{(}\bigcup_{k=n_0}^{n}\{2b\bar{S}_{k}\geq
b^{2}\bar{V}^{2}_{k}+b^{2}\ep(\bar{V}^{2}_{n}-\bar
{V}^{2}_{k})+x^{2}+\varepsilon\} \Biggr{)}\\[-2pt]
&&\qquad{}-\pr\Biggl{(}\bigcup_{k=n_0}^{n}\{2b\bar{S}_{k}\geq b^{2}\bar{V}^{2}_{k}+
x^{2}+\varepsilon, (\bar{V}^{2}_{n}-\bar{V}^{2}_{k})-\ep
(\bar{V}^{2}_{n}-\bar{V}^{2}_{k})\geq\varepsilon B^{2}_{n}/x^{2}\}
\Biggr{)}\nn\\[-2pt]
&&\quad=:I_{5}-I_{6}.\nn
\end{eqnarray}
To complete the proof of Proposition~\ref{lemma-2}, we only need to
show the following lemma.
\begin{lemma}\label{lemma3} Under the conditions of Proposition \ref
{lemma-2}, we have
\begin{eqnarray}
\label{nd1}I_{5}&\geq&2 \bigl(1-\Phi(x) \bigr) \bigl(1-C_{\delta, \tau}
(\varepsilon^{-2}\Delta_{n,x}+\varepsilon) \bigr),\\[-2pt]
\label{nd2}I_{6}&\leq& C_{\tau,\delta}x^{-2}\mathrm{e}^{-x^{2}/2}.
\end{eqnarray}
\end{lemma}
\begin{pf*}{Proof of (\ref{nd1})}
We have
\begin{eqnarray}\label{a10}
I_{5}&\geq& \pr\Biggl{(}\bigcup_{k=1}^{n}\{2b\bar{S}_{k}\geq
b^{2}\bar{V}^{2}_{k}+b^{2}\ep(\bar{V}^{2}_{n}-\bar
{V}^{2}_{k})+x^{2}+\varepsilon\} \Biggr{)}\no\\[-2pt]
& &{}- \pr\Biggl{(}\bigcup_{k=1}^{n_0}\{
2b\bar{S}_{k}\geq
b^{2}\bar{V}^{2}_{k}+b^{2}\ep(\bar{V}^{2}_{n}-\bar
{V}^{2}_{k})+x^{2}+\varepsilon\} \Biggr{)}\\[-2pt]
&=:&I_{5,1}-I_{5,2}.\no
\end{eqnarray}
Write
\begin{eqnarray*}
\xi_{j}&=&2b\bar{X}_{j}-b^{2}\bar{X}^{2}_{j}+b^{2}\ep\bar
{X}^{2}_{j},\\[2pt]
\F_{j}&=&\Biggl\{\sum_{k=1}^{j}\xi_{k}\geq
y\Biggr\},\qquad \mbox{where } y=2x^{2}+\varepsilon.
\end{eqnarray*}
As in the proof of (\ref{a18}), there exists a constant $c_0$ such
that for all $0\le k\le n-1$,
\begin{eqnarray*}
\pr\Biggl(\sum_{j=k+1}^{n}(\xi_{j}-\ep\xi_{j}) \geq c_{0}\varepsilon
\Biggr)\leq1/2.
\end{eqnarray*}
This, together with the independence of $\xi_j$, yields that
\begin{eqnarray*}
I_{5,1}&=&\pr(\F_{1})+\sum_{k=2}^{n}\pr(\F^{c}_{1},\ldots, \F
^{c}_{k-1},\F_{k})\\[2pt]
&\geq& \pr(\F_{1},y\leq\xi_{1}\leq
y+4\varepsilon)+\sum_{k=2}^{n}\pr\Biggl(\F^{c}_{1},\ldots, \F
^{c}_{k-1},\F
_{k}, y\leq\sum_{j=1}^{k}\xi_{j}\leq y+4\varepsilon\Biggr)\cr
&\geq&2\pr\Biggl(\F_{1}, y\leq\xi_{1}\leq
y+4\varepsilon, \sum_{j=2}^{n}(\xi_{j}-\ep\xi_{j}) \geq
c_{0}\varepsilon\Biggr)\\[2pt]
& &{}+2\sum_{k=2}^{n}\pr\Biggl(\F^{c}_{1},\ldots, \F^{c}_{k-1},\F_{k},
y\leq
\sum_{j=1}^{k}\xi_{j}\leq y+4\varepsilon,\sum_{j=k+1}^{n}(\xi
_{j}-\ep\xi
_{j}) \geq c_{0}\varepsilon\Biggr)\\[2pt]
&\geq&2\pr\Biggl{(}\sum_{k=1}^{n}(\xi_{k}-\ep\xi_{k})\geq
y+(c_{0}+4)\varepsilon+D_{n} \Biggr{)},
\end{eqnarray*}
where $D_{n}=\sum_{j=1}^{n}|\ep\xi_j|$.
Similarly to the proofs of (\ref{a19})--(\ref{a21}), it follows from
Lemma~\ref{le3} with $\eta_j=\xi_j-\ep\xi_j$ that
\begin{eqnarray}\label{a8}
I_{5,1} &\geq& 2\pr\Biggl{(}\sum_{k=1}^{n}(\xi_{k}-\ep\xi_{k})\geq
y+(c_{0}+4)\varepsilon+D_{n} \Biggr{)}\no
\\[-8pt]
\\[-8pt]
&\geq&2
\bigl(1-\Phi(x)\bigr ) \bigl(1-A (\varepsilon+\varepsilon^{-2}\Delta_{n,x}) \bigr).\no
\end{eqnarray}
On the other hand, similar to the proofs of (\ref{a19}) and (\ref
{a21}), we have
\begin{eqnarray}\label{a9}
I_{5,2}&\leq& 2\pr\Biggl{(}\sum_{j=1}^{n_0}\xi_{j}\geq
2x^{2}+(1-c_{0})\varepsilon-D_{n} \Biggr{)}\no\\
&\leq&Cx^{-1}\exp\biggl{(}-\frac{x^{2}}{2}-\frac{x^{2}\sum
_{j=n_0+1}^{n}\ep
X^{2}_{j}}{2B_{n}} \biggr{)}\\
&\leq&Cx^{-2}\mathrm{e}^{-x^{2}/2}.\no
\end{eqnarray}
This, together with (\ref{a8}), implies (\ref{nd1}).\vadjust{\goodbreak}
\end{pf*}

\begin{pf*}{Proof of (\ref{nd2})}
Define $k'_{0}=1$, and
$k'_{i}=k'_{i-1}+1$ if $\ep X^{2}_{k'_{i-1}+1}> \varepsilon^2
B^{2}_{n}/x^{6}$, and otherwise
\[
k'_{i}=\max\Biggl{\{}k\leq n\dvt \sum_{j=k'_{i-1}+1}^{k}\ep
X^{2}_{j}\leq
\frac{\varepsilon^{2}B^{2}_{n}}{x^{6}} \Biggr{\}}+1.
\]
Let $m$ satisfy $k'_{m-1}<n\leq k'_{m}$ and define
\[
k_{i}=k'_{i} \qquad\mbox{for }i<m,\quad \mbox{and}\quad k_{m}=n.
\]
Because $\sum_{j=k_{i-1}+1}^{k_{i}}\ep X^{2}_{j}> \varepsilon
^{2}B^{2}_{n}/ x^{6}$ for $i<m$, we have
\[
B^{2}_{n}\geq\sum_{i=1}^{m-1}\sum_{j=k_{i-1}+1}^{k_{i}}\ep
X^{2}_{j}>(m-1)\varepsilon^{2}B^{2}_{n}/x^{6},
\]
which implies that $m\leq\varepsilon^{-2}x^{6}+1$. Furthermore,
suppose that $i_{0}$ satisfies $k_{i_{0}-1}<n_{0}\leq k_{i_{0}}$,
where $n_0$ is defined as in (\ref{n0a}). Set
\begin{eqnarray*}
\breve{X}_{k} &= & X_{k}I\{|X_{k}|\leq16^{-1}\varepsilon B_{n}/x^{3}\}, \\
\hat{X}_{k} & = & X_{k}I\{16^{-1}\varepsilon B_{n}/x^{3}<|X_{k}|\leq
\varepsilon B_{n}/x\}, \qquad\hat{Z}_{k_{i}}=\sum
_{k=k_{i-1}+1}^{k_{i}-1}|\hat
{X}_{k}|.
\end{eqnarray*}
Note that $2b|\bar{X}_{k}|\leq2\varepsilon$. Simple calculations
show that
\begin{eqnarray}\label{df0}
I_6
&\leq& \sum_{i=i_{0}}^{m}\pr\Biggl{(}
\bigcup_{k=k_{i-1}+1}^{k_{i}}\{2b\bar{S}_{k}\geq b^{2}\bar{V}^{2}_{k}+
x^{2}+\varepsilon, \nn\\
&&\hphantom{\sum_{i=i_{0}}^{m}\pr\Biggl{(}
\bigcup_{k=k_{i-1}+1}^{k_{i}}\{}
(\bar{V}^{2}_{n}-\bar{V}^{2}_{k})-\ep(\bar{V}^{2}_{n}-\bar
{V}^{2}_{k})\geq\varepsilon B^{2}_{n}/x^{2}\} \Biggr{)}\no\\
&\leq&
\sum_{i=i_{0}}^{m}\pr\Biggl(\bigcup_{k=k_{i-1}+1}^{k_{i}-1}\{2b\bar
{S}_{k}\geq b^{2}\bar{V}^{2}_{k}+ x^{2}-\varepsilon,\cr
&& \hphantom{\sum_{i=i_{0}}^{m}\pr\Biggl{(}\bigcup_{k=k_{i-1}+1}^{k_{i}-1}\{}(\bar{V}^{2}_{n}-\bar{V}^{2}_{k_{i-1}+1})-\ep(\bar
{V}^{2}_{n}-\bar
{V}^{2}_{k_{i-1}+1})\geq2^{-1}\varepsilon
B^{2}_{n}/x^{2}\} \Biggr)\no
\\[-8pt]
\\[-8pt]
&\leq& \sum_{i=i_{0}}^{m}\pr\bigl(2b(\bar{S}_{k_{i-1}}+\hat
{Z}_{k_{i}}+\ep\hat{Z}_{k_{i}})\geq b^{2}\bar{V}^{2}_{k_i}+
x^{2}-2\varepsilon,\no\\[-2pt]
&&\hphantom{\sum_{i=i_{0}}^{m}\pr{(}}(\bar{V}^{2}_{n}-\bar{V}^{2}_{k_{i-1}+1})-\ep
(\bar{V}^{2}_{n}-\bar{V}^{2}_{k_{i-1}+1})\geq2^{-1}\varepsilon
B^{2}_{n}/x^{2} \bigr)\no\\[-2pt]
&&{} + \sum_{i=i_{0}}^{m}\pr\Biggl( \max_{k_{i-1}+1\leq
j\leq k_{i}-1}2b\sum_{k=k_{i-1}+1}^{j}(\breve{X}_{k}-\ep\breve
{X}_{k})\geq\varepsilon\Biggr)\no\\[-2pt]
&=:&I_{6,1}+I_{6,2}. \no
\end{eqnarray}
Noting that
\[
\sigma_{ni}^2:=
\sum_{k=k_{i-1}+1}^{k_{i}-1}\ep\breve{X}^{2}_{k}\leq
\frac{\varepsilon^{2} B^{2}_{n}}{x^{6}} \quad\mbox{and} \quad|\breve
{X}_{k}|\leq
16^{-1}\varepsilon
B_{n}/x^{3},
\]
it follows from $m\leq\varepsilon^{-2}x^{6}+1$ and L\'evy's
inequality that
with $t=2bx^2/\varepsilon$
\begin{eqnarray}\label{df1}
I_{6,2} &\le&
\sum_{i=i_{0}}^{m}\pr\Biggl{(}
\sum_{k=k_{i-1}+1}^{k_i-1}(\breve{X}_{k}-\ep\breve{X}_{k})\geq
\varepsilon/(2b)-\sqrt2 \sigma_{ni} \Biggr{)} \no\\[-2pt]
&\le& \sum_{i=i_{0}}^{m} \mathrm{e}^{-t(\varepsilon/(2b)-\sqrt2 \sigma_{ni})}
\prod_{k=k_{i-1}+1}^{k_i-1}E\mathrm{e}^{t(\breve{X}_{k}-\ep\breve{X}_{k})}
\nn
\\[-8pt]
\\[-8pt]
&\leq& A\mathrm{e}^{-x^2} \sum_{i=i_{0}}^{m}\exp\{At^2 \sigma_{ni}^2 \} \no
\\[-2pt]
&\leq& 2A_1(\varepsilon^{-2}x^{6}+1) \mathrm{e}^{-x^{2}}\leq
C_{\tau,\delta}x^{-2}\mathrm{e}^{-x^{2}/2},\no
\end{eqnarray}
where we used the fact that $\varepsilon\ge\gamma x^{-1/2}$.

To estimate $I_{6,1}$, let $t=24 \varepsilon^{-1} x^{2}B_n^{-2}\log x$.
Note that
\[
2b\ep\hat{Z}_{k_{i}}\leq\frac{32x^{4}}{\varepsilon
B^{2}_{n}}\sum_{k=k_{i-1}+1}^{k_{i}-1}\ep X^{2}_{k}\leq
\frac{32\varepsilon}{x^{2}}\leq8\varepsilon.
\]
Similar to the estimate for $I_3$ in (\ref{ad12}), we obtain
\begin{eqnarray}\label{last1}
I_{6,1} &\leq&
\sum_{i=i_{0}}^{m}\pr\bigl{(}2b(\bar{S}_{k_{i-1}}+\hat{Z}_{k_{i}}-\ep
\hat{Z}_{k_{i}})\geq b^{2}\bar{V}^{2}_{k_{i-1}}+
x^{2}-18\varepsilon,\no\\[-2pt]
&&\hphantom{\sum_{i=i_{0}}^{m}\pr\bigl{(}}
(\bar{V}^{2}_{n}-\bar{V}^{2}_{k_{i-1}+1})-\ep
(\bar{V}^{2}_{n}-\bar{V}^{2}_{k_{i-1}+1})\geq2^{-1}\varepsilon
B^{2}_{n}/x^{2} \bigr{)}\no\\[-2pt]
&\leq&
\sum_{i=i_{0}}^{m}\pr\bigl{(}b(\bar{S}_{k_{i-1}}-\ep
\bar{S}_{k_{i-1}}+\hat{Z}_{k_{i}}-\ep
\hat{Z}_{k_{i}})-b^{2}\bar{V}^{2}_{k_{i-1}}/2 \nn\\[-2pt]
&&\hphantom{\sum_{i=i_{0}}^{m}\pr\bigl{(}}{} +t(\bar{V}^{2}_{n}-\bar{V}^{2}_{k_{i-1}+1})-t\ep
(\bar{V}^{2}_{n}-\bar{V}^{2}_{k_{i-1}+1})
\geq x^{2}/2+12\log
x-9\varepsilon\bigr{)}\no\\[-2pt]
&\leq&
Ax^{-12}\mathrm{e}^{-x^{2}/2}\sum_{i=i_{0}}^{m}
\Biggl\{ \prod_{j=1}^{k_{i-1}} \ep
\mathrm{e}^{b(\bar{X}_{j}-\ep
\bar{X}_{j})-2^{-1}b^{2}\bar{X}^{2}_{j}} \nn
\\[-8pt]
\\[-8pt]
&&\phantom{Ax^{-12}\mathrm{e}^{-x^{2}/2}\sum_{i=i_{0}}^{m}
\Biggl\{}{} \times\prod_{j=k_{i-1}+1}^{k_{i}-1}\ep
\mathrm{e}^{b(|\hat{X}_{j}|-\ep|\hat{X}_{j}|)+t(\bar{X}^{2}_{j}-\ep
\bar{X}^{2}_{j})}
\times
\prod_{j=k_{i}}^{n}\ep \mathrm{e}^{t(\bar{X}^{2}_{j}-\ep\bar{X}^{2}_{j})} \Biggr\}
\no\hspace*{25pt}
\\[-2pt]
&\leq& Ax^{-12}\mathrm{e}^{-x^{2}/2}
\sum_{i=i_{0}}^{m}\exp\biggl{(}A\Delta_{n,x}+A\frac{x^{2}\sum
_{j=k_{i-1}+1}^{k_{i}-1}\ep
X^{2}_{j}}{B^{2}_{n}}\mathrm{e}^{24\varepsilon\log x}\no\\[-2pt]
&&\hphantom{Ax^{-12}\mathrm{e}^{-x^{2}/2}
\sum_{i=i_{0}}^{m}\exp\biggl{(}}{} +A\frac{x^{3}\sum_{j=k_{i-1}+1}^{k_{i}-1}\ep
|\bar{X}_{j}|^3}{B^{3}_{n}}\mathrm{e}^{24\varepsilon\log
x}\varepsilon^{-1}\log x\no\\[-2pt]
&&\hphantom{Ax^{-12}\mathrm{e}^{-x^{2}/2}
\sum_{i=i_{0}}^{m}\exp\biggl{(}}{}
+A\frac{x^{4}\sum_{j=k_{i-1}+1}^{k_{i}-1}\ep
\bar{X}^{4}_{j}}{B^{4}_{n}}\mathrm{e}^{24\varepsilon\log
x}(\varepsilon^{-1}\log x)^{2}\no\\[-2pt]
&&\hphantom{Ax^{-12}\mathrm{e}^{-x^{2}/2}
\sum_{i=i_{0}}^{m}\exp\biggl{(}}{}
+A\frac{x^{4}\sum_{j=k_{i}}^{n}\ep
\bar{X}^{4}_{j}}{B^{4}_{n}}\mathrm{e}^{24\varepsilon\log
x}(\varepsilon^{-1}\log x)^{2} \biggr{)}.\no
\end{eqnarray}
Recall, by the definition of $k_i$,
\begin{eqnarray*}
\frac{x^2}{B_n^2}\sum_{j=k_{i-1}+1}^{k_{i}-1}\ep{X}_{j}^2 &\le
&\varepsilon^2 x^{-4}, \\[-2pt]
\frac{x^4}{B_n^4} \sum_{j=k_{i-1}+1}^{k_{i}-1}\ep
\bar{X}^{4}_{j}&\leq& \frac{\varepsilon x^3}{B_n^3}
\sum_{j=k_{i-1}+1}^{k_{i}-1}\ep|\bar{X}_{j}|^3 \le \frac
{\varepsilon^2 x^2}{B_n^2} \sum_{j=k_{i-1}+1}^{k_{i}-1}\ep
{X}_{j}^2 \le\varepsilon^4 x^{-4}.
\end{eqnarray*}
On the other hand, we have
\[
\varepsilon^{-1}x^{24\varepsilon}\le\gamma^{-1}\min\{x^{\delta
/10+\delta/3}, x^{3/2}\}\le
\gamma^{-1}x^{\min\{\delta/2, 3/2\}},
\]
and by (\ref{n0a})--(\ref{a2}) and the
inequality $\sum_{j=n_0+1}^{n}\ep X_{j}^{2}\leq192x^{-2}B^{2}_{n}\log x$,
for all $i\ge i_0$,
\begin{eqnarray*}
\frac{x^4}{B_n^4}\sum_{j=k_{i}}^{n}\ep\bar{X}^{4}_{j}&\leq&
(\varepsilon B_n/x) \sum_{j=n_0+1}^{n}\ep|\bar{X}_{j}|^{3}\\[-2pt]
& \leq&(\varepsilon\tau B_n^2/x^{2+\delta}) \sum_{j=n_0+1}^{n}\ep
X_{j}^{2} \leq C_{\tau,\delta}\varepsilon x^{-\delta} \log x.
\end{eqnarray*}
Substituting these estimates into (\ref{last1}) gives
\be
I_{6,1} &\leq& C_{\delta, \tau} (\varepsilon^{-2}x^{6}+1)
x^{-12}\mathrm{e}^{-x^{2}/2}\exp\(C_{\tau,\delta}x^{-1/2}\log^2 x+C_{\tau
,\delta
}x^{-\delta/2}\log^3 x\) \no
\\[-9pt]
\\[-9pt]
&\le& C_{\delta, \tau} x^{-2}\mathrm{e}^{-x^{2}/2}.\no
\ee
This proves (\ref{nd2}), which also completes the proof of Lemma~\ref{lemma3}.\vadjust{\goodbreak}
\end{pf*}

\section*{Acknowledgements}
Weidong Liu's research partially supported by a foundation  of national excellent doctoral
dissertation of PR China. Qi-Man Shao's research partially supported by Hong Kong
RGC CERG 602608 and 603710. Qiying Wang's research partially supported by an ARC discovery grant.
We thank the Associate Editor and the referee for their helpful
comments, which led to a
significant improvement of the presentation of the paper.
%
%

\printhistory
\end{document}